\documentclass[12pt]{amsart}

\usepackage{amsmath, amssymb}
\usepackage[unicode]{hyperref}
\usepackage{tikz}
\usepackage{pgfplots}
\usepackage{filecontents}

\newtheorem{Thm}{Theorem}

\begin{document}

\title{Integers without large prime factors: from Ramanujan to de Bruijn}

\author{Pieter Moree}
\address{Max-Planck-Institut f\"ur Mathematik,Vivatsgasse 7, D-53111 Bonn, Germany}
%\curraddr{}
\email{moree@mpim-bonn.mpg.de}
%\thanks{}

\date{}


\dedicatory{In memoriam: Nicolaas Govert (`Dick') de Bruijn (1918-2012) }

%\dedicatory{In memoriam: Nicolaas Govert (`Dick') de Bruijn (1918-2012) }
%{\def\thefootnote{}
%\footnote{{\it Mathematics Subject Classification (2000)}.
%11T22, 11B83}}

\maketitle

\begin{abstract}
\noindent A small survey of work done
on estimating the number of integers without large prime factors up to around 1950 is
provided. Around 1950 N.G. de Bruijn published results that dramatically advanced the subject and
started a new era in this topic.
\end{abstract}
\section{Introduction}
Let $P(n)$ denote the largest prime divisor
of $n$. 
We set $P(1)=1$.
A number $n$ is said to be {\it $y$-friable}\footnote{In the older
literature one usually finds $y$-smooth. Friable
is an adjective meaning easily crumbled or broken.} if $P(n)\le y$. We let $S(x,y)$ denote the set of integers $1\le n\le x$
such that $P(n)\le y$. The cardinality of $S(x,y)$ is denoted by $\Psi(x,y)$. We write 
$y=x^{1/u}$, that is $u=\log x/\log y$.\\
\indent Fix $u>0$. In 1930, Dickman \cite{Dickman} proved that
\begin{equation}
\label{dikkertje}
\lim_{x\rightarrow \infty}{\Psi(x,x^{1/u})\over x}=\rho(u),
\end{equation}
with
$$
\rho(u)=\rho(N)-\int_N^u {\rho(v-1)\over v}dv,~~(N<u\le N+1,~N=1,2,3,\ldots),
$$
and $\rho(u)=1$ for $0<u\le 1$ (see Figure 1). 
It is left to the reader to show that we have
\begin{equation}
\label{defie}
\rho(u)=
\begin{cases}
1 & \text{for $0\le u\le 1$};\\
{1\over u}\int_{0}^1 \rho(u-t)dt & \text{for $u>1$.}
\end{cases}
\end{equation}
The function $\rho(u)$ in the literature is either called the {\it Dickman function} 
or the {\it Dickman-de Bruijn function}.\\ 
\indent In this note I will briefly discuss the work done on friable integers up to the papers of 
de Bruijn \cite{most,B3} that appeared around 1950 and dramatically advanced the subject. A lot of the early work
was carried out by number theorists from India.\\
\indent De Bruijn \cite{most} improved on (\ref{dikkertje}) by establishing a result that together with
the best currently known estimate for the prime counting function (due to I.M. Vinogradov and Korobov in 1958) yields
the following result.
\begin{Thm}
\label{vijf}
The estimate 
\begin{equation}
\label{bruni}
\Psi(x,y)=x\rho(u)\Big\{1+O\footnote{The reader not familiar with the Landau-Bachmann O-notation we refer to wikipedia or any
introductory text on analytic number theory, e.g., Tenenbaum \cite{T}. Instead of $\log \log x$ we sometimes write $\log_2 x$, instead of
$(\log x)^A$, $\log^A x$.}\Big({\log(u+1)\over \log y}\Big)\Big\},
\end{equation}
holds for $1\le u\le \log^{3/5-\epsilon}y,{~that~is,~}y>\exp(\log^{5/8+\epsilon}x)$.
\end{Thm}
De Bruijn's most important tool in his proof of this result is the {\it Buchstab equation} \cite{Buchstab},
\begin{equation}
\label{boek}
\Psi(x,y)=\Psi(x,z)-\sum_{y<p\le z}\Psi({x\over p},p),
\end{equation}
where $1\le y<z\le x$. The Buchstab equation is easily proved on noting that the number of integers $n\le x$ with $P(n)=p$
equals $\Psi(x/p,p)$. Given a good estimate for $\Psi(x,y)$ for $u\le h$, it allows one to obtain a good
estimate for $u\le h+1$.\\ 
\indent De Bruijn \cite{B3} complemented Theorem \ref{vijf} by an asymptotic estimate for $\rho(u)$.
That result has as a corollary that, for $u\ge 3$,
\begin{equation}
\label{eerstdick}
\rho(u)=\exp\Big\{-u\Big\{\log u+\log_2u-1+{\log_2 u-1\over \log u}+O\Big(\big({\log_2u\over \log u}\big)^2\Big)\Big\}\Big\},
\end{equation}
which will suffice for our purposes. 
Note that (\ref{eerstdick}) implies that, as $u\rightarrow \infty$,
$$\rho(u)={1\over u^{u+o(u)}},~~\rho(u)=\Big({e+o(1)\over u\log u}\Big)^u,$$
formulas that suffice for most purposes and are easier to remember.
For a more detailed description of this and other work
of de Bruijn in analytic number theory, we refer to Moree \cite{engels}.

\section{Results on $\rho(u)$}
Note that $\rho(u)>0$, for if not, then because of the continuity of $\rho(u)$ there is a smallest zero $u_0>1$ 
and then substituting $u_0$ in (\ref{defie}) we easily arrive at a contradiction.
Note that for $u>1$ we have 
\begin{equation}
\label{deriva}
\rho'(u)=-{\rho(u-1)\over u}
\end{equation} 
It follows that $\rho(u)=1-\log u$ for $1\leq u\leq 2$. 
For $2\leq u \leq  3$, $\rho(u)$ can be expressed in terms of the dilogarithm. However, with increasing
$u$ one has to resort to estimating $\rho(u)$ or finding a numerical approximation.\\
\indent Since $\rho(u)>0$ we see from (\ref{deriva}) that $\rho(u)$ is strictly decreasing for $u>1$.
{}From this and (\ref{defie}) we then find that $u\rho(u)\le \rho(u-1)$, which on using
induction leads to
$\rho(u)\le 1/[u]!$ for $u\ge 0$. 
It follows that $\rho(u)$ quickly tends to
zero as $u$ tends to infinity.\\ 
\indent Ramaswami \cite{rma2} proved that
$$\rho(u)>{C\over u 4^u \Gamma(u)^2},~u\ge 1,$$
for a suitable constant $C$, with $\Gamma$ the Gamma function. 
By Stirling's formula we have $\log \Gamma (u)\sim u\log u$ and
hence the latter inequality is for $u$ large enough improved on by the following inequality
due to Buchstab \cite{Buchstab}:
\begin{equation}
\label{buchen}
\rho(u)>\exp\Big\{-u\Big\{\log u+\log_2u+6{\log_2 u\over \log u}\Big\}\Big\},~~(u\ge 6).
\end{equation}
Note that on its turn de Bruijn's result (\ref{eerstdick}) considerably improves on the latter inequality.
\begin{figure}
\caption{The Dickman-de Bruijn function $\rho(u)$}
\medskip
\begin{tikzpicture}
  \begin{axis}[
    width=\linewidth,
    height=5cm,
    axis x line=middle,
    axis y line=middle,
    xlabel=$u$, 
    xmin=-0.1, xmax=4.2,
    ymin=-0.05, ymax=1.1
    ]
    \addplot[color=blue] file {rho.dat};
    \node at (axis cs:2,0.5) {$\rho(u)$};
  \end{axis}
\end{tikzpicture}
\end{figure}
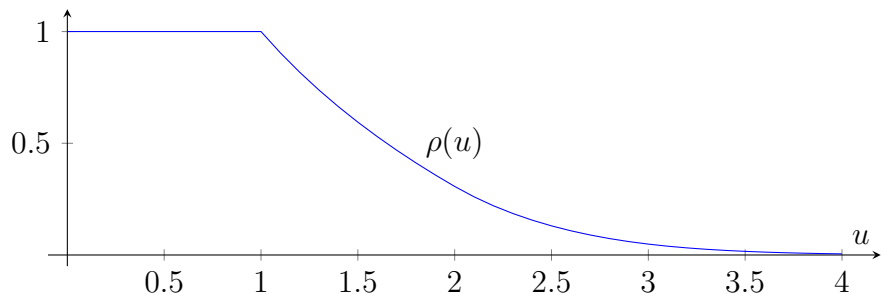

\section{S. Ramanujan (1887-1920) and the friables}
In his first letter (January 16th, 1913) to Hardy (see, e.g. \cite{BR}), one of the most famous letters in all of mathematics, Ramanujan
claims that
\begin{equation}
\label{oudebrief}
\Psi(n,3)={1\over 2}{\log (2n)\log (3n)\over \log 2\log 3}.
\end{equation}
The formula is of course intended as an approximation, and there is no evidence to show how
accurate Ramanujan supposed it to be. Hardy \cite[pp. 69-81]{Hardybook} in his lectures on Ramanujan's work 
gave an account of an interesting analysis that can be made to hang upon the above assertion.
I return to this result in the section on the $\Psi(x,y)$ work of Pillai.\\
\indent In the so-called Lost Notebook \cite{Ramanujan} we find at the second half of page 337: 

$\phi(x)$ {\sl is the no. of nos of the form} 
$$ 
2^{a_2} \cdot 3^{a_3} \cdot 5^{a_5}\cdots p^{a_p} \qquad p \le x^{\epsilon}
$$
{\sl not exceeding $x$}.  
$$ 
\tfrac 12\le \epsilon \le 1, \qquad \phi(x) \sim x \Big\{ 1 -\int_\epsilon^1 \frac{d\lambda_0}{\lambda_0}\Big\}
$$
$$ 
\tfrac 13\le \epsilon \le \tfrac 12, \qquad \phi(x) \sim x \Big\{ 1-\int_{\epsilon}^1 \frac{d\lambda_0}{\lambda_0} + \int_{\epsilon}^{\frac 12} \frac{d\lambda_1}{\lambda_1} \int_{\lambda_1}^{1-\lambda_1}\frac{d\lambda_0}{\lambda_0}\Big\}
$$
\begin{eqnarray*} 
\tfrac 14 \le \epsilon \le \tfrac {1}{3}, \qquad \phi(x)\sim x \Big \{ 1&&- \int_{\epsilon}^1 \frac{d\lambda_0}{\lambda_0} + \int_{\epsilon}^{\frac 12} \frac{d\lambda_1}{\lambda_1} \int_{\lambda_1}^{1-\lambda_1}\frac{d\lambda_0}{\lambda_0} \\
&&- \int_{\epsilon}^{\frac 13} \frac{d\lambda_2}{\lambda_2} \int_{\lambda_2}^{\frac{1-\lambda_2}{2}} 
\frac{d\lambda_1}{\lambda_1} \int_{\lambda_1}^{1-\lambda_1} \frac{d\lambda_0}{\lambda_0}\Big\} \\
\end{eqnarray*}
\begin{eqnarray*}
\tfrac 15 \le \epsilon \le \tfrac 14, \qquad \phi(x) \sim x &&\Big\{ 1- \int_{\epsilon}^1 \frac{d\lambda_0}{\lambda_0} + \int_{\epsilon}^{\frac 12} \frac{d\lambda_1}{\lambda_1} \int_{\lambda_1}^{1-\lambda_1}\frac{d\lambda_0}{\lambda_0} \\
&&- \int_{\epsilon}^{\frac 13} \frac{d\lambda_2}{\lambda_2} \int_{\lambda_2}^{\frac{1-\lambda_2}{2}} 
\frac{d\lambda_1}{\lambda_1} \int_{\lambda_1}^{1-\lambda_1} \frac{d\lambda_0}{\lambda_0}\\
&&+\int_{\epsilon}^{\frac 14}
\frac{d\lambda_3}{\lambda_3} \int_{\lambda_3}^{\frac{1-\lambda_3}{3}} \frac{d\lambda_2}{\lambda_2} 
\int_{\lambda_2}^{\frac{1-\lambda_2}{2}} \frac{d\lambda_1}{\lambda_1} \int_{\lambda_1}^{1-\lambda_1} \frac{d\lambda_0}{\lambda_0}\Big\}\\
\end{eqnarray*}
{\sl and so on.} \\

\indent In the book by Andrews and Berndt \cite[\S 8.2]{AB} it is
shown that Ramanujan's assertion is equivalent with (\ref{dikkertje}) with
$$ 
\rho(u) = \sum_{k=0}^{\infty} \frac{(-1)^{k}}{k!} I_k(u), 
$$ 
where 
$$ 
I_k(u) = \int_{{t_1, \ldots t_k \ge 1}\atop{ t_1+\ldots+t_k \le u}} \frac{dt_1}{t_1} \cdots 
\frac{dt_k}{t_k}. 
$$ 
This is one of many examples where Ramanujan reached with his hand from his grave to snatch a theorem, in this
case from Dickman who was at least 10 years later than Ramanujan, cf. Berndt \cite{Bruce}.
Chowla and Vijayaragahavan \cite{CV} seemed to have been the first to rigorously prove
(\ref{dikkertje}) with $\rho(u)$ expressed as a sum of iterated integrals (cf. the section
on Buchstab). The asymptotic behaviour of the integrals $I_k(u)$ has been studied by
Soundararajan \cite{sound}.\\
\indent Ramanujan's claim reminds me of
the following result of 
Chamayou \cite{cham}: If $x_1,x_2,x_3,\cdots$ are independent random variables uniformly distributed in $(0,1)$, and 
$u_n=x_1+x_1x_2+\ldots+x_1x_2\cdots x_n$, then $u_n$ converges in probability to a limit $u_{\infty}$ and 
$u_{\infty}$ has a probability distribution with density function $\rho(t)e^{-\gamma}$, where 
$\gamma$ denotes Euler's constant. 

\section{I.M. Vinogradov (1891-1983) and the friables}
The first to have an application for $\Psi(x,y)$ estimates seems to have been 
Ivan Matveyevich Vinogradov \cite{vino} in 1927.
Let $k\ge 2$ be a prescribed integer and $p\equiv 1({\rm mod~}k)$ a prime. The $k$-th powers in $(\mathbb Z/p\mathbb Z)^*$ form
 a subgroup of order $(p-1)/k$ and so the existence follows of 
$g_1(p,k)$, the least $k$-th power non-residue modulo a prime $p$.
Suppose that $y<g_1(p,k)$, then $S(x,y)$ consists of $k$-th power residues only. It follows that
$$\Psi(x,y)\le \# \{n\le x~:~n\equiv a^k({\rm mod~}p){\rm~for~some~}a\}.$$  
The idea is now to use good estimates for the quantities on both side of the inequality sign in order to deduce an upper bound
for $g_1(p,k)$.\\
\indent Vinogradov \cite{vino} showed that $\Psi(x,x^{1/u})\ge \delta(u)x$ for $x\ge 1$, $u>0$, where 
$\delta(u)$ depends only on $u$ and is positive. He applied this to show that if $m\ge 8$, $k>m^m$, and $p\equiv 1({\rm mod~}k)$
is sufficiently large, then 
\begin{equation}
\label{expo}
g_1(p,k)<p^{1/m}.
\end{equation}
See Norton \cite{Norton} for a historical account of the problem
of determining $g_1(p,k)$ and original results.

\section{K. Dickman (1861-1947) and the friables}
Karl Dickman was active in the Swedish insurance business in the end of the 19th century and the beginning of the 20th century.
Probably, he studied mathematics in the 1880's at Stockholm University, where the legendary Mittag-Leffler
was professor\footnote{I have this information from Lars Holst.}.\\
\indent As already mentioned Dickman proved (\ref{dikkertje}) and in the
same paper\footnote{Several sources falsely claim that Dickman wrote only one mathematical
paper. He also wrote \cite{second}.} gave an heuristic argument to the effect that
\begin{equation}
\label{karl}
\lim_{x\rightarrow \infty}{1\over x}\sum_{2\leq n\leq x}{\log P(n)\over \log n}=\int_0^{\infty}{\rho(u)\over (1+u)^2}du.\end{equation}
Denote the integral above by $\lambda$. Dickman argued that $\lambda \approx 0.62433$.
Mitchell \cite{Mitchell} in 1968 computed that $\lambda=0.62432 99885 4\dots$.
The interpretation of Dickman's heuristic is that
for an average integer with $m$ digits, its greatest prime factor has about $\lambda m$ digits.
The constant $\lambda$ is now known as
the {\it Golomb-Dickman constant}, as it arose independently in research of Golomb and others involving the largest
cycle in a random permutation.\\
\indent De Bruijn \cite{most} in 1951 was the first to prove (\ref{karl}). He did
this using his $\Lambda(x,y)$-function, an approximation
of $\Psi(x,y)$, that he introduced in the same paper.

\section{S.S. Pillai (1901-1950) and the friables}
Subbayya Sivasankaranarayana Pillai (1901-1950) was a number theorist who worked on problems in classical
number theory (Diophantine equations, Waring's problem, etc.). 
Indeed, he clearly was very much inspired by the work of Ramanujan.
He tragically died in a plane crash near Cairo
while on his way to the International Congress of
Mathematicians (ICM) 1950, which was held at Harvard University.\\
\indent Pillai wrote
two manuscripts on friable integers, \cite{Pillaiun1, Pillaiun2}, of which \cite{Pillaiun1} was accepted for
publication in the Journal of the London Mathematical Society, but did not appear in print. Also \cite{Pillaiun2} was
never published in a journal.\\
\indent In \cite{Pillaiun1}, see also \cite[pp. 481-483]{Pillai2}, Pillai investigates $\Psi(x,y)$ for $y$ fixed. Let $p_1,p_2,\ldots,p_k$ denote all the 
different primes $\le y$. Notice
that $\Psi(x,y)$ equals the cardinality of the set
$$\{(e_1,\ldots,e_k)\in \mathbb Z^k: e_i\ge 0,~\sum_{i=1}^k e_i\log p_i\le x\}.$$
Thus $\Psi(x,y)$ equals the number of lattice points in a $k$-dimensional tetrahedron with sides of
length $\log x/\log 2,\ldots, \log x/\log p_k$. This tetrahedron has volume
$${1\over k!}\prod_{p\le y}\big({\log x\over \log p}\big).$$
Pillai shows that 
$$\Psi(x,y)={1\over k!}\prod_{p\le y}\big({\log x\over \log p}\big)\Big(1+(1+o(1)){k\log(p_1p_2\ldots p_k)\over 2\log x}\Big).$$
If $\rho_1,\ldots,\rho_k$ are positive real numbers and $\rho_1/\rho_2$ is irrational, then the same estimate with
$\log p_i$ replaced by $\rho_i$ holds for 
$$\{(e_1,\ldots,e_k)\in \mathbb Z^k: e_i\ge 0,~\sum_{i=1}^k e_i\rho_i\le x\}.$$
This was proved by Specht \cite{Specht} (after whom the Specht modules are named),
see also Beukers \cite{Frits}.
A much sharper result than that of Pillai/Specht was obtained in 1969 by 
Ennola \cite{Ennola} (see also Norton \cite[pp. 24-26]{Norton}). In this result Bernoulli
numbers make their appearance.\\
\indent Note that Pillai's result implies that
\begin{equation}
\label{soep}
\Psi(x,3)={1\over 2}{\log (2x)\log (3x)\over \log 2\log 3}+o(\log x),
\end{equation}
and that the estimate
$$\Psi(x,3)={\log^2 x\over 2\log 2\log 3}+o(\log x)$$
is false. Thus Ramanujan's estimate (\ref{oudebrief}) is more precise than 
the trivial estimate $\log^2 x/(2\log 2\log 3)$. Hardy \cite[\S 5.13]{Hardybook}
showed that the error term $o(\log x)$ in
(\ref{soep}) can be
replaced by $o(\log x/\log_2 x)$. In the proof of this he uses a result of Pillai \cite{PH1}, see
also \cite[pp. 53-61]{Pillai1}, saying that given $0<\delta<1$,
one has $|2^x-3^y|>2^{(1-\delta)x}$ for all integers $x$ and $y$ with $x>x_0(\delta)$ sufficiently large.\\
\indent In \cite{Pillaiun2}, see also \cite[pp. 515-517]{Pillai2}, Pillai claims that, for $u\ge 6$, ${B/u}<\rho(u)<{A/u}$,
with $0<B<A$ constants.
He proves this result by induction assuming a certain estimate for $\rho(6)$ holds.
However, this estimate for $\rho(6)$ does not hold. Indeed, the claim contradicts (\ref{eerstdick}) and
is false.\\
\indent Since Pillai reported on his work on the friables at conferences in India and 
stated open problems there, his
influence on the early development of the topic was considerable. E.g., one of the questions he raised was whether
$\Psi(x,x^{1/u})=O(x^{1/u})$ uniformly for $u\le (\log x)/\log 2$. This question was answered in the affirmative
by Ramaswami \cite{rma2}.

\section{R.A. Rankin (1915-2001) and the friables}
In his work on the gaps between consecutive primes Robert Alexander Rankin \cite{Rankin} in 1938 introduced a simple
idea to estimate $\Psi(x,y)$ which turns out to be remarkably effective and can be used in similar situations. 
This idea is now called `Rankin's method' or `Rankin's trick'.
Starting
point is the observation that for any $\sigma>0$
\begin{equation}
\label{rankinne}
\Psi(x,y)\le \sum_{n\in S(x,y)}({x\over n})^{\sigma}\le x^{\sigma}\sum_{P(n)\le y}{1\over n^{\sigma}}=x^{\sigma}\zeta(\sigma,y),
\end{equation}
where 
$$\zeta(s,y)=\prod_{p\le y}(1-p^{-s})^{-1},$$
is the partial Euler product up to $y$ for the Riemann zeta function $\zeta(s)$. Recall that, for $\Re s>1$,
$$\zeta(s)=\sum_{n=1}^{\infty}{1\over n^s}=\prod_{p}{1\over 1-p^{-s}}.$$
By making an appropriate choice for $\sigma$ and estimating $\zeta(\sigma,y)$ using analytic prime number
theory, a good upper bound for $\Psi(x,y)$ can be found. E.g., the choice
$\sigma=1-1/(2 \log y)$ leads to
$$\zeta(\sigma,y)\ll \exp \Big\{\sum_{p\leq y}{1\over p^{\sigma}}\Big\}\leq \exp\Big\{
\sum_{p\leq y}{1\over p}+O\Big((1-\sigma)\sum_{p\leq y}{\log p\over p}\Big)\Big\}\ll \log y,$$
which gives rise to 
\begin{equation}
\label{rancune}
\Psi(x,y)\ll x{\rm e}^{-u/2}\log y.
\end{equation}

\section{A.A. Bukhshtab (1905-1990) and the friables}
Aleksandr Adol'fovich Bukhshtab (Buchstab in the German spelling) most important contribution is
the equation (\ref{boek}) now named after him. A generalization of it plays an important role in
sieve theory. Buchstab \cite{Buchstab} in 1949 proved (\ref{dikkertje})
and gave both Dickman's differential-difference equation as well as the result
\begin{equation}
\label{akin}
\rho(u)=1+\sum_{n=1}^N (-1)^n\int_n^u \int_{n-1}^{t_1-1}\int_{n-2}^{t_2-1}\cdots \int_1^{t_{n-1}-1}{dt_ndt_{n-1}\cdots dt_1\over
t_1t_2\cdots t_n},
\end{equation}
for $N\le u\le N+1$ and $N\ge 1$ an integer, simplifying Chowla and Vijayaragahavan's expression (they erroneously
omitted the term $n=N$). Further, Buchstab
established inequality (\ref{buchen}) and applied his results to show that the exponent in 
Vinogradov's result (\ref{expo}) can be roughly divided by two.

\section{V. Ramaswami and the friables}
V. Ramaswami\footnote{He worked at Andhra University until his death in 1961. I will be grateful for further bibliographical information.} \cite{rma1} showed that
$$\Psi(x,x^{1/u})=\rho(u)x+O_U({x\over \log x}),$$
for $x>1$, $1<u\le U$, and remarked that the error term is best possible. He sharpened this  result
in \cite{rma2} and showed there that, for $u>2$,
\begin{equation}
\label{vrijdag}
\Psi(x,x^{1/u})=\rho(u)x+\sigma(u){x\over \log x}+O({x\over \log^{3/2}x}),
\end{equation}
with $\sigma(u)$ defined similarly to $\rho(u)$. Indeed, it turns out that $$\sigma(u)=(1-\gamma)\rho(u-1),$$ but this
was not noticed by Ramaswami. In \cite{rma3} Ramaswami generalized his results to $B_l(m,x,y)$ which counts
the number of integers $n\le x$ with $P(n)\le y$ and $n\equiv l({\rm mod~}m)$
\footnote{Buchstab \cite{Buchstab} was the first to investigate $B_l(m,x,y)$.}.
Norton \cite[pp. 12-13]{Norton} points out 
some deficits of this paper and gives a reproof \cite[\S 4]{Norton} of Ramaswami's result on $B_l(m,x,x^{1/u})$ generalizing
(\ref{vrijdag}).\\
\indent {}From de Bruijn's paper \cite[Eqs. (5.3), (4.6)]{most} one easily derives the following generalization
of Ramaswami's results\footnote{The notation $O_m$ indicates that the implied constant might depend on $m$.}:
\begin{Thm}
\label{twee}
Let $m\ge 0$, $x>1$, and suppose $m+1<u<\sqrt{\log x}$. Then
$$\Psi(x,y)=x\sum_{r=0}^m a_r{\rho^{(r)}(u)\over \log^r y}+{O_m}\big({x\over \log^{m+1}y}\big),$$
with $\rho^{(r)}(u)$ the $r$-th derivative of $\rho(u)$ and $a_0,a_1,\ldots$ are the coefficients
in the power series expansion
$${z\over 1+z}\zeta(1+z)=a_0+a_1z+a_2z^2+\ldots,$$ with $|z|<1$.
\end{Thm}
It is well-known (see, e.g., Briggs and Chowla \cite{BC}) that around $s=1$ the Riemann zeta function has 
the Laurent series expansion
$$\zeta(s)={1\over s-1}+\sum_{k=0}^{\infty}{(-1)^k\over k!}\gamma_k(s-1)^k,$$
with 
$\gamma_k$ the $k$-th Stieltjes constant and with $\gamma_0=\gamma$ Euler's constant.
Using this we find that $a_0=1$ and $a_1=\gamma-1$. Thus Theorem \ref{twee} yields (\ref{vrijdag}) with
$\sigma(u)=(1-\gamma)\rho(u-1)$ for the range $2<u<\sqrt{\log x}$. For $u>\sqrt{\log x}$ 
the estimate (\ref{vrijdag}) in view of (\ref{eerstdick}) reduces to
$$\Psi(x,x^{1/u})\ll {x\log^{-3/2}x},$$ which easily follows from (\ref{rancune}).

\section{S. Chowla (1907-1995) and the friables}
The two most prominent number theorists in the period following Ramanujan were S.S. Pillai and
Sarvadaman Chowla. They kept in contact through an intense correspondence \cite{sutha}. Chowla in his long career
published hunderds of reseach papers.\\
\indent  Chowla and Vijayaragahavan \cite{CV} expressed $\rho(u)$ as an iterated integral and gave a formula akin
to (\ref{akin}). De Bruijn \cite{B1} established some results implying that $\Psi(x,\log^h x)=O(x^{1-1/h+\epsilon})$ for $h>2$.
An easier reproof of the latter result was given by Chowla and Briggs \cite{CB}. 

\section{Summary}
It seems that the first person to look at friable integers was Ramanujan, starting with his first letter
to Hardy (1913), also Ramanujan seems to have been the first person to arrive at the Dickman-de Bruijn function
$\rho(u)$. Pillai generalized some of Ramanujan's work and spoke about it on conferences in India, which
likely induced a small group of Indian number theorists to work on friable integers. 
Elsewhere in the same period (1930-1950) only
incidental work was done on the topic. Around 1950 N.G. de Bruijn published his ground-breaking 
papers \cite{most, B3}. Soon afterward the
Indian number theorists stopped publishing on friable integers.\\
\indent It should also be said that the work on friable
integers up to 1950 seems to contain more mistakes than more recent work.
Norton \cite{Norton} points out and corrects many of these mistakes.\\

\noindent {\tt Further reading}. 
As a first introduction to friable numbers I can highly recommend Granville's 2008 survey \cite{Granville}. 
It has a special
emphasis on friable numbers and their role in algorithms in computational number theory. 
Mathematically more demanding is
the 1993 survey by Hildebrand and Tenenbaum \cite{HT1}. 
Chapter III.5 in Tenenbaum's book \cite{T}
deals with $\rho(u)$ and approximations to $\Psi(x,y)$ by the saddle point method.\\

\noindent {\tt Acknowledgement}. 
I thank R. Thangadurai for helpful correspondence on S.S. Pillai and the friables and providing me with a PDF file of Pillai's collected works. 
B.C. Berndt kindly sent me a copy of \cite{AB}. 
K.K. Norton provided helpful comments on an earlier version. His research monograph \cite{Norton}, which
is the most extensive source available on the early history of friable integer counting, was quite helpful to me. In \cite{Norton}, by the way,
new results (at the time) on $g_1(p,k)$ and $\Psi_m(x,y)$, the number of $y$-friable integers $1\le n\le x$ coprime with
$m$ are established.
Figure 1 was kindly created for me by Jon Sorenson and Alex Weisse (head of the MPIM computer group).

\end{document}